\documentclass[12pt]{amsart}

\usepackage{amsmath,amssymb}

\newtheorem{thm}{Theorem}
\newtheorem{prop}[thm]{Proposition}
\newtheorem*{rem}{Remark}

\newcommand{\Z}{{\mathbb Z}}
\newcommand{\Q}{{\mathbb Q}}
\newcommand{\R}{{\mathbb R}}

\newcommand{\prob}{\operatorname{Prob}}
\newcommand{\trace}{\operatorname{Tr}}
\newcommand{\real}{\operatorname{Re}}
\newcommand{\Frob}{\operatorname{Frob}}
\newcommand{\fp}{{\mathbb F}_p}
\newcommand{\fq}{{\mathbb F}_q}
\newcommand{\fpext}{{\mathbb F}_{p^2}}
\newcommand{\fpk}{{\mathbb F}_{p^k}}
\newcommand{\PP}{{\mathbb P}}
\newcommand{\FF}{{\mathcal F}}
\newcommand{\OO}{{\mathcal O}}
\newcommand{\mm}{{\mathfrak m}}

\newcommand {\E} {\mathbb{E}}

\title[Gaussian point count statistics for families of curves]
{Gaussian point count statistics for families of curves over a fixed finite field}

\author{P\"ar Kurlberg}
\urladdr{www.math.kth.se/\~{ }kurlberg}
\address{Department of Mathematics, Royal Institute of Technology,
 SE-100 44 Stockholm, Sweden}
 \email{kurlberg@math.kth.se}

\author{Igor Wigman}
\urladdr{www.math.kth.se/\~{ }wigman}
\address{Department of Mathematics, Royal Institute of Technology,
 SE-100 44 Stockholm, Sweden}
 \email{wigman@kth.se}

\thanks{P.K. was partially supported by grants from
  the Knut and Alice Wallenberg foundation, the Royal Swedish Academy
  of Sciences, and the Swedish Research Council.  I.W.  was supported
  by grant KAW 2005.0098 from the Knut and Alice Wallenberg
  Foundation.}

\begin{document}
\begin{abstract}
  We produce a collection of families of curves, whose point count
  statistics over $\fp$ becomes Gaussian for $p$ {\em fixed}.  In
  particular, the {\em average} number of $\fp$-points on curves in
  these families tends to infinity.
\end{abstract}
\date{June 23, 2010}
\maketitle


\section{Introduction}
\label{sec:introduction}


The purpose of this note is to exhibit a collection of families of smooth
curves whose normalized limiting point count statistics, over a {\em
  fixed} finite field $\fp$, have a Gaussian distribution.
Given a finite family $\FF_i$ of smooth curves defined over $\fp$, let
$
M_{i} := \frac{1}{|\FF_i|} \sum_{C \in \FF_{i}} |C(\fp)|
$
be the average number of $\fp$-points on the curves $C \in \FF_{i}$
and let
$
V_{i} := \frac{1}{|\FF_i|} \sum_{C \in \FF_{i}} \left( |C(\fp)|  - M_{i} \right)^{2}
$
be the variance of the fluctuations in these point counts.
Here, and in what follows, $C(\fp)$ will denote the set of
$\fp$-points on $C$, and $|C(\fp)|$ its cardinality.
We can now formulate the main result of this note.
\begin{thm}
\label{thm:main}
There exists a sequence of families $\{\FF_i\}_{i=1}^{\infty}$ of
smooth curves defined over $\fp$ with the following properties:
$|\FF_i|, M_{i}, V_{i}$ all tend to infinity,
and, for all compact intervals $I$,
$$
\frac{1}{|\FF_i|} \left | \left\{  C \in \FF_{i}:   \frac{|C(\fp)|-M_i}{V_i^{1/2}} \in I
  \right\} \right |
= \frac{1}{\sqrt{2\pi}}\int\limits_{I} e^{-x^2/2}dx + o(1),
$$
as $i \to \infty$.
\end{thm}


To obtain such a sequence
we intersect projective
surfaces $X_{i} \subset \PP^{n_{i}}$, chosen so that $|X_{i}(\fp)|$
tends to infinity as $i \to
\infty$, with families of large degree hypersurfaces. More precisely,
let $S_{i}(d)$ be the set of degree $d$ homogeneous polynomials in
$n_{i}+1$ variables with coefficients in $\fp$.  For $f \in S_{i}(d)$,
let $H_{f} \subset \PP^{n_{i}}$ be the hypersurface defined by the zero
set of $f$.  Intersecting $X_i$ with $H_f$ we generically obtain a (possibly
singular) curve
$
C_{i}(f) := X_{i} \cap H_{f}.
$
Letting
$$
\tilde{S}_{i}(d) := \{ f \in S_i(d) : X_{i} \cap H_{f} \text{ is smooth} \},
$$
we obtain a family of {\em smooth} curves
$$
\FF_{i}(d) := \{ C_{i}(f) : f \in \tilde{S}_{i}(d) \}.
$$
%

%


Our main technical result (in essence a slightly more explicit version
of Poonen's \cite[Theorem~1.2]{poonen-finite-field-bertini}) asserts
that the distribution of point counts for curves in this family, for
$d$ large, is binomial --- we can think of it as the number of
successes when an unfair coin is tossed $|X_i(\fp)|$ times.
\begin{prop}
\label{prop:main-technical}
Let $X_i \subset \PP^{n_i}$ be a smooth projective surface defined
over $\fp$, and let $\FF_i(d)$ be defined as above.  Then, as $d \to \infty$,
\begin{multline}
\label{eq:C pnts <-> binom}
\frac{| \{  C \in \FF_i(d):   |C(\fp)| = s   \}|}
{|\FF_i(d)|}
\\=
\binom{|X_i(\fp)|}{s}
\left(\frac{p+1}{p^{2}+p+1}\right)^{s}
\left(1-\frac{p+1}{p^{2}+p+1}\right)^{|X_i(\fp)|-s}
\cdot (1+o(1))
\end{multline}
{\em uniformly} for $0 \le s \le |X_i(\fp)|$. In particular, the average
point count of a curve $C \in \FF_i(d)$ equals $|X_{i}(\fp)| \cdot
\frac{(p+1)}{p^{2}+p+1} \cdot (1+o(1))$ as $d \to \infty$.
\end{prop}

Given Proposition~\ref{prop:main-technical}, we can easily obtain
Theorem~\ref{thm:main}
by the central limit theorem type argument for coin flip models (cf.
Section~\ref{sec:main-thm-proof}) {\em provided} we
can find a sequence of surfaces $\{X_{i}\}_{i \geq1}$ such that
$|X_{i}(\fp)| \to \infty$.
Any such sequence
suffices, but for concreteness we will use Ihara's construction
\cite{ihara81-some-remarks} (independently discovered by Tsfasman,
Vl{\u{a}}du{\c{t}}, and Zink
\cite{tsfasman-vladut-zink-curves-with-many-points}), of families of
{\em curves} with many points over $\fpext$, and a (Weil) restriction
of scalars argument will then produce a projectively embedded
$\fp$-surface having many points --- see
Section~\ref{sec:surfaces-with-many} for more details.

{\em Remark:} Since the relation $|C(\fp)| = p + 1 -
\trace(\Frob|H_c^1(\overline{C}, \overline{\Q}_l))$ holds between the
number of $\fp$-points on a curve $C$ and the trace of Frobenius (see
e.g \cite[Ch. 11]{IwaniecKowalsk-analytic-NT-book}), we have in fact
exhibited families of smooth curves for which the average of the trace
of Frobenius exhibits a strong negative bias.  (Note that the average
trace of Frobenius should be zero according to random matrix theory
predictions; cf. Section~\ref{sec:gaussian-point-count}.)

\newpage
\subsection{Background and discussion}
\subsubsection{Gaussian point count statistics in other models}
\label{sec:gaussian-point-count}

The number of $\fp$-points on a smooth curve $C$ of genus $g$, defined
over a finite field $\fp$,
can be written as
$$
|C(\fp)| = p+1 - \sum_{i=1}^{2g} \alpha_i
$$
where $\alpha_i \in \bar{\Q}$ are the eigenvalues of the Frobenious
action on a certain cohomology group.  By Weil's proof of the Riemann
hypothesis for curves, $|\alpha_i| = \sqrt{p}$ for $1 \leq i \leq 2g$.
Hence $|C(\fp)| = p+1 - p^{1/2} \cdot \trace(U_{C})$, where $U_{C} \in
U({2g})$ is a unitary $2g \times 2g$ matrix, unique up to conjugacy.
If we let $C$ range over a family $\FF(\fp)$ of smooth curves defined
over $\fp$ (e.g., the family of hyperelliptic curves $\FF = \{ C_{t}
\}_{t \in \fp : f(t) \neq 0}$ where $C_{t} = \{ (x,y): y^{2} =
f(x)(x-t)\}$ and $f \in \fp[x]$ is a square free polynomial of degree
$2g$) it is natural to study the distribution of the fluctations of
the points counts by looking at the normalized fluctuations, i.e., the
quantity $(|C(\fp)|
- p - 1)/\sqrt{p} = -\trace(U_{C})$.

By Deligne's equidistribution theorem, the distribution of the
conjugacy classes of the $U_{C}$'s, as $p \to \infty$, are given by
random matrix theory when the family has ``large monodromy''.  For
example, Katz and Sarnak \cite{KatSar99} has shown that for the family
of hyperelliptics given above, the limiting distribution on the
$U_{C}$-conjugacy classes, as $p \to \infty$, is given by the Haar
measure on $USp(2g)$, the group of unitary symplectic $2g\times2g$
matrices.  On the other hand, in the limit $g \to \infty$, Diaconis
and Shahshahani has shown \cite{DiaSha94} that the limiting distribution of
$\{\trace(U)\}_{U \in USp(2g)}$, as $g \to \infty$, is a {\em
  Gaussian} with mean zero and variance one\footnote{It is rather
  remarkable that even though $\trace(U)$ is a sum of $2g$ complex
  numbers on the unit circle, the variance does not scale as
  $\sqrt{g}$, but is in fact identically equal to one.}.  Thus, if a
collection of families of curves have large monodromies, then the
(normalized) point count fluctuations are Gaussian in the double limit
$\lim_{g \to \infty} \lim_{p \to \infty}$.

For $p$ fixed and $g\to \infty$ it is less clear what to expect
regarding the distribution of the $U_{C}$-conjugacy classes and their
traces.  For instance, random matrix theory (RMT) is clearly not an
appropriate model since the inequality $ 0 \leq |C(\fp)| = p+1 -
\sqrt{p} \cdot \trace(U_C)$ does not hold for all $U_{C} \in USp(2g)$
for $g$ sufficiently large\footnote{However, it is worth noting that
  certain statistics of the eigenvalues are consistent with RMT, e.g.,
  the fluctuations of the number of eigenvalues in random short
  intervals (cf.  \cite{faifman-rudnick08}).  Moreover, in
  \cite{Rudnick10-traces-of-frobenius} Rudnick found that the
  one-level density, a local statistics, was in agreement with RMT
  (even though averages of traces of small powers was not.)}.
Similarly, if the curves in the families can be embedded into
$\PP^n(\fp)$ for $n$ fixed, the bounds $0 \leq |C(\fp)| = p+1 -
\sqrt{p} \cdot \trace(U_C) \leq |\PP^{n}(\fp)|$ rules out a Gaussian.
In fact the normalized distribution cannot even have continuous
support since $|C_{f}(\fp)|$ is integer valued and the variance is
bounded.


To get some insight into the large genus limit while keeping the
ground field fixed, Kurlberg and Rudnick studied
(cf. \cite{ap-distribution}) families of hyperelliptic curves of the
form $C_{f} : y^{2} = f(x)$, where $f$ ranges over monic square free
polynomials of degree $d$.  They found that the fluctuations in this
family, as $d\to\infty$, has the same distribution as $\sum_{i=1}^{p}
x_{i}$, where $x_{1}, \ldots, x_{p}$ are independent random variables
taking the values $0,-1,1$ with probabilities $1/(p+1), 1/2(1+p^{-1}),
1/2(1+p^ {-1})$, respectively.  Moreover, the moments for the
normalized point count distribution were shown to be Gaussian as long
as both $p,d \to \infty$, i.e., any joint limit rather than letting $p$
tend to infinity first (as described above.)

In \cite{bucur-etal-cyclic-covers-imrn}, Bucur, David, Feigon and
Lal\'in generalized this to cyclic $l$-fold covers of $\PP^{1}$; here
the distribution of the fluctuations are given by $2 \real(
\sum_{i=1}^{p} x_{i})$ where $x_{1}, \ldots, x_{p}$ are independent
random variables taking the values $0,\exp(2 \pi i/l), \ldots, \exp(2
\pi i (l-1)/l)$ with probabilities $2/(p+2), p/(3(p+2)), \ldots,
p/(3(p+2))$, respectively.  Further, in
\cite{bucur-etal-plane-curves}, they studied the family of smooth
curves in the projective plane cut out by degree $d$ homogenous
polynomials, and found that the distribution of the point count
statistics for this family, as $d\to \infty$, is the same as that of
$\sum_{i=1}^{p^{2}+p+1} x_{i}$, where $x_{1}, \ldots, x_{p^{2}+p+1}$
are independent random variables taking the values $0,1$ with
probabilities $p^{2}/(p^{2}+p+1), (p+1)/(p^{2}+p+1)$, respectively.
%
In both \cite{bucur-etal-cyclic-covers-imrn,bucur-etal-plane-curves},
Gaussian moments were obtained when $p$ tends to infinity; in the
first case along any limit $p,d \to \infty$, whereas the assumption $d
\gg p^{1+\epsilon}$ is needed in the second case.

We finally note that in \cite{Lar-gauss-preprint} (unpublished),
Larsen obtained Gaussian moments for a smooth family of hyperelliptic
curves of the form $Y^{2} = \prod_{i=1}^{n}(X-a_{i})$, where $a_{1},
\ldots, a_{n} \in \fq$, ranges over distinct elements.

\subsubsection{Remarks on vanishing probabilities}
Given a point $P \in X_i(\fp)$, the probability of a polynomial
$f \in S_{i}(d)$ vanishing at $P$ is $1/p$, so one might
expect that the average number of $\fp$-points on $X_{i} \cap H_{f}$
should equal $|X_i(\fp)|/p$.
However, as we have seen, this prediction is not quite correct --- by
conditioning on $f$ so that $X_i \cap H_f$ is smooth,
the probability of $f$ vanishing at $P$ turns out to be slightly smaller than
expected, and is given by $\frac{p+1}{p^{2}+p+1}$ (rather than by
$1/p$). A similar phenomenon was already observed in \cite{bucur-etal-plane-curves} for
smooth plane curves: the probability of point $P \in \PP^2(\fp)$
belonging to a smooth curve given by the zero set of a random
homogenous polynomial (of large degree) is
$(p+1)/(p^{2}+p+1)$, rather than $1/p$.

\subsubsection{Remarks on families of curves with many points}

In order to obtain Gaussian (normalized) point count statistics, it is
essential that there is no a priori upper bound on the number of
$\fp$-points on the curves; in particular, families of plane curves,
families with bounded genus, or families with bounded gonality cannot
be used.
In fact, something even stronger is needed: since point counts are
integer valued, the variance must grow to infinity for the normalized
distribution to have continuous support.  Thus, since the limiting
distribution must be symmetric around the mean, together with the fact
that the number of points on a curve is non-negative, the {\em
  average} number of $\fp$-points of a curve in the family also must
tend to infinity.

A natural candidate for families of curves with unbounded point counts
over $\fp$ is $M_g(\fp)$, the set of isomorphism classes of genus $g$
curves.
However, {\em if} these point counts
can be
modeled by random matrix theory (e.g., as in the case of $g$ fixed and
$p \to \infty$ as shown in \cite{KatSar99}), the average number of
points 
would be $p+1$ since the average trace
of Frobenius equals zero by random matrix theory
predictions; hence it is unclear whether using $M_g(\fp)$,
$g\to\infty$ as a collection of families would work.  (Also see
\cite{brock-granville-excess, katz-quadratic-excess} for explicit
results on average point counts, as $p \to \infty$, for curves in
various families.)


Another possibility to avoid a priori upper bounds might be to
consider smooth curves given by intersecting $n-1$ generic hypersurfaces in
$\PP^n$.
However, the average number of points on curves in this family turns
out to be bounded; Bucur and Kedlaya recently
\cite{bucur-kedlaya-smooth-intersection} showed that it is slightly
{\em less} than $p+1$.  In particular, the average trace of Frobenius
for curves in this family is {\em not} equal to zero.

\subsubsection{Acknowledgments:}
The authors would like to thank Alina Bucur,
Chantal David, Nicholas Katz, and Ze\'ev Rudnick for helpful
discussions.

\section{Proof of Theorem~\ref{thm:main}}
\label{sec:main-thm-proof}
\subsection{Normal distribution from coin flip model}

We first recall some facts about the binomial distribution.
For $i\ge 1$, let $B_{i}$ be the binomial random variable counting
the successful tosses among
\begin{equation}
\label{eq:ni def}
n_{i} := |X_i(\fp)|
\end{equation}
tosses of an unfair
coin, whose probability of success equals
\begin{equation}
\label{eq:r prob def}
r:=\frac{p+1}{p^2+p+1}.
\end{equation}
Let $M_{i}'$ and $V_{i}'$ be the expected value
and the variance of the number $B_{i}$ of successful tosses
respectively. That is,
$
M_i' := \E[B_{i}] =n_{i} \cdot r
$
and
$
V_i' := Var(B_{i}) =  n_{i} \cdot r(1-r).
$
It is a classical result in probability (see e.g. ~\cite{feller},
chapter 7) that
$B_{i}$, suitably normalized, tends to the standard Gaussian,
provided that $n_{i}\rightarrow\infty$ (note that the probability of a success
stays constant).
Namely,
$
\frac{B_{i}-M_{i}'}{\sqrt{V_{i}'}} \rightarrow \mathcal{N}(0,1)
$
in the sense that for any compact interval $I\subseteq\R$,
\begin{equation}
\label{eq:CLT binom}
\prob\left( \frac{B_{i}-M_{i}'}{\sqrt{V_{i}'}} \in I  \right)
\rightarrow \frac{1}{\sqrt{2\pi}}\int\limits_{I} e^{-x^2/2}dx.
\end{equation}

We are now in  position to prove Theorem \ref{thm:main}.

\begin{proof}[Proof of Theorem \ref{thm:main} assuming Proposition \ref{prop:main-technical}]

Recall the definitions of $n_{i}$ and $r$ in Equations \eqref{eq:ni def} and \eqref{eq:r prob def} respectively.
For the coin flip model the probability of precisely $k$ successes is given by
$$
\prob\left( B_{i} = k \right) =
\binom{n_{i}}{k}
r^{k}
\left(1-r\right)^{n_{i}-k}
$$
whereas, by Proposition~\ref{prop:main-technical},
$\prob\left( |C(\fp)| = k  \right) $ for our family of curves
equals
$$
\binom{n_{i}}{k}
r^{k}
\left(1-r\right)^{n_{i}-k}
\cdot (1+o(1))
$$
with constant involved in the `$o$'-notation being uniform in $k$. In
particular, we find that
$$
V_{i} =
\frac{1}{|\FF_i|} \sum_{C \in \FF_{i}} \left( |C(\fp)|  - M_{i}
\right)^{2}
=
V_{i}' \cdot (1+o(1)).
$$
Further,
\begin{multline*}
  \prob\left(\frac{|C(\fp)| - M_{i}}{V_{i}} \in I \right)
 =
\sum_{\substack{k \in \Z  \\ k  \in M_{i} + I \cdot V_{i} \cap [0,n_{i}]}}
\binom{n_{i}}{k}
r^{k}
\left(1-r\right)^{n_{i}-k} \cdot (1+o(1))
\end{multline*}
so by comparing with the corresponding sum in the coin flip model 
and using the classical \eqref{eq:CLT binom}, together with
$$
\sum_{k=0}^{n_{i}}
\binom{n_{i}}{k}
r^{k}
\left(1-r\right)^{n_{i}-k}
=1
$$
to control the error term, we obtain the statement of Theorem
\ref{thm:main}.

\end{proof}

\subsection{Surfaces with many points}
\label{sec:surfaces-with-many}

A crucial assumption in Theorem~\ref{thm:main} is the existence of a
sequence of surfaces $\{X_{i}\}$ whose point counts over $\fp$ tends
to infinity; here we give a concrete example of such a sequence.
We begin by recalling the  construction of modular curves with many
$\fpext$-points given in
\cite{tsfasman-vladut-zink-curves-with-many-points}.
\begin{thm}
If $l$ is a prime greater than $p$, there exists a smooth complete
curve\footnote{$X_0(l)$ is the projective smooth model of the
  modular curve $Y_{0}(l)$, which parametrizes pairs $(E,C)$ where $E$
  is an elliptic curve and $C \subset E$ is a cyclic subgroup of order $l$.}
$X_{0}(l)$, defined over $\fp$,
having at least $(p-1)(l+1)/12$ points over $\fpext$.
\end{thm}
If $C$ is curve of genus $g>1$, and $C$ is not hyperelliptic, there
exists a canonical embedding $\phi : C \to \PP^{g-1}$, whose image is
a curve of degree $2g-2$ inside $\PP^{g-1}$
(cf. \cite[Ch.~4.5]{MR0463157}).  Now, the genus of $X_{0}(l)$ equals
$g_{l} = [l/12]$, and since a hyperellipic curve has at most $2(p+1)$
points, $X_{0}(l)$ clearly cannot be hyperelliptic for $l$
sufficiently large.  Thus $X_{0}(l)$ can be canonically embedded into
$\PP^{g_{l}-1}$ for $l$ large.

Letting $l_{i}$ denote the $i$-th prime we now define a sequence of
surfaces $\{ X_{i} \}$ by letting $X_{i}$ be the
$\fp$-surface obtained by Weil restriction of scalars, from $\fpext$
to $\fp$, of $X_0(l_{i})$.
As remarked above,
$X_{0}(l_{i})/\fpext$ is canonically embedded into
$\PP^{g_{l_{i}}-1}_{\fpext}$, and it is known (e.g., see Lemma~7.5 in
Ch. I.7.2 of
\cite{Handbook-elliptic-hyperelliptic-curve-cryptography})
that $\PP^{n}_{\fpext}$ can be projectively embedded in
$\PP^{(n+1)^{2}-1}_{\fp}$, hence each $X_{i}$ is a projective surface
with $|X_{i}(\fp)| = |X_{0}(l_{i})(\fpext)| \geq (p-1)(l_{i}+1)/12$ for
$i$ sufficiently large.

\section{Proof of Proposition~\ref{prop:main-technical}}

We begin with some notation: given a variety $X$ defined over $\fq$,
let $\zeta_{X}(s)$ denote the zeta function attached to $X$, i.e.,
$$
\zeta_{X}(s) := \prod_{
\substack{
P \in X \\ \text{$P$ closed}
}
} (1-p^{-\deg(P)s})^{-1}.
$$

For $C_{i}(f) = X_{i} \cap H_{f}$ to be smooth, $H_{f}$ must intersect
$X_{i}$ transversally at all points $P \in X_{i}(\overline{\fp})$ such that
$f(P)=0$.  I.e., if $f(P)=0$ for $P \in X(\fpk)$ and we write
$f|_{X_{i}}$ in local coordinates $x,y$ (say with $P$ corresponding to
$x=y=0$), we must have $f|_{X_{i}} = ax+by+(\text{higher order
  terms})$, where $a,b \in \fpk$ {\em and} $(a,b) \neq (0,0)$.  In
other words, $f|_{X_{i}}$ must not have quadratic order of vanishing
at any point $P \in X_{i}(\fpk)$ (for any $k$), or equivalently, the
image of $f$ in $\OO_{X,P}/\mm_{X,P}^{2}$ must be nonzero
for all $P \in X_{i}(\fpk)$ and all $k$.

Since $X$ is a smooth surface, $|\OO_{X,P}/\mm_{X,P}^{2}| = p^{3k}$
for all $P \in X_{i}(\fpk)$, so the ``probability of smoothness at
$P$'' equals $1-p^{-3k}$; if these conditions are sufficently
independent as $P$ varies, the probability of $C_{i}(f)$ being smooth
should be given by $1/\zeta_{X_i}(3)$.  Using a sieving argument,
Poonen showed that this heuristic indeed gives the correct answer when
$d\to \infty$ and $f$ ranges over elements in $S_{i}(d)$.  Further,
his ``Finite field Bertini with Taylor coefficients'' \cite[Theorem
1.2]{poonen-finite-field-bertini} allows for controlling the behaviour
of $f$ in the neighborhood of a finite number of points, best
formulated in terms of schemes.  We shall need the following slightly
more explicit version:
\begin{thm}
\label{thm:bertini-with-taylor}
Let $X$ be a quasiprojective subscheme of $\PP^n$ over $\fq$, and let
$S(d)$ be the set of degree $d$ homogenuous polynomials in $n+1$
variables.   Let $Z$
be a finite subscheme of $\PP^n$, and assume that $U := X - (Z \cap
X)$ is smooth of dimension $2$.  Fix a subset $T \subset
H^0(Z,\OO_Z)$.  Given $f \in S(d)$, let $f|_Z$ be the element of
$H^0(Z,\OO_Z)$ that on each connected component $Z_k$ equals the
restriction of $x_j^{-d}f$ to $Z_k$, where $j=j(k)$ is the smallest $j
\in {0,1,\ldots, n}$ such that the coordinate $x_j$ is invertible on
$Z_k$. Then as $d\rightarrow\infty$
\begin{multline*}
|\{ f \in S(d) : H_f \cap U \text{ is a smooth curve, and } f|_Z \in
  T \}|
\\
= \frac{|S(d)|}{\zeta_U(3)} \cdot \frac{|T|}{|H^0(Z, \OO_Z)|}
\cdot (1+o_{Z}(1)),
\end{multline*}
\end{thm}
\begin{proof}
We will closely follow the closed point sieve of
\cite[Section~2]{poonen-finite-field-bertini}.
The case $T = \emptyset$ is trivial, so we may assume that $|T| \geq 1$.
Let $U_{r}$ be the set of closed points of $U$ of degree $<r$, and
define $U_{ \geq r}$ similarly.
By the proof of \cite[Lemma 2.2]{poonen-finite-field-bertini}, for $r$
fixed and $d$ sufficiently large,
\begin{multline*}
|\{ f \in S(d) : H_{f} \cap U \text{ is smooth at all
  $P \in U_{r}$ and $f|_{Z} \in T$} \}|
\\ =
 \frac{|S(d)|  \cdot |T|}{|H^{0}(Z,\OO_Z)|} \cdot
 \prod_{P \in U_{<r}} (1-p^{-3 \deg(P)})
\end{multline*}
which in turn equals
$$
\frac{|S(d)|  \cdot |T|}{ |H^{0}(Z,\OO_Z)|} \cdot
\frac{ (1 + O(p^{-3r}))}{\zeta_{U}(3)}
$$
Letting $r$ grow with $d$ in a suitable fashion, we obtained the
claimed main term; to conclude the proof it is enough to bound the
number $f \in S(d)$ for which smoothness of $H_{f} \cap U$ is violated
at some point $P$ of degree larger than $r$.

{\em Medium degree points:} By the proof of
\cite[Lemma~2.4]{poonen-finite-field-bertini},
$$
|\{ f \in S(d) : H_{f} \cap U \text{ is not smooth at some $P \in U_{\geq
    r} \cap U_{< d/3}$}  \}|
\ll_{U} |S(d)| \cdot p^{-r}.
$$

{\em Large degree points:} By the proof of
\cite[Lemma~2.6]{poonen-finite-field-bertini} (in particular, see
Claim~1 and Claim~2), there exists $\tau \in \Z^+$, only depending on
$U$, such that
\begin{multline*}
|\{ f \in S(d) : H_{f} \cap U \text{ is not smooth at some $P \in
  U_{\geq d/3}$} \}|
\\=
|S(d)| \cdot O( d p^{-(d-\tau)/p} + d^{2} p^{-\min([d/p]+1,d/3)}   )
\end{multline*}
which, for $d$ sufficiently large, is
$$
\ll
|S(d)| p^{-d/(p+2)}
$$
Now, taking $r = [d/4]$, we find that the total contribution from
medium and large degree primes is
$
\ll |S(d)| p^{-d/5}
$
and the result follows by taking $d$ sufficiently large so that
$
p^{-d/(p+3)} = o \left(\frac{1}{|H^{0}(Z,\OO_Z)| \zeta_U(3)  }   \right)
$.

\end{proof}

\begin{rem}
In particular, taking  $Z$ to be empty, we obtain
\begin{multline}
\label{eq:number-of-smooth-guys}
|\tilde{S}(d)| =
|\{ f \in S(d) : H_f \cap X \text{ is a smooth curve} \}|
\\
= \frac{|S(d)|}{\zeta_X(3)}
\cdot (1+o(1)),
\end{multline}
which can be interpreted as saying that the probability of $H_{f} \cap
X$ being smooth is $1/\zeta_X(3)+o(1)$ for $X$ a surface.
\end{rem}
\subsection{Applications to point counting}
To apply Theorem~\ref{thm:bertini-with-taylor}, we define a finite
subscheme $Z \subset X$ as follows: let $X(\fp) =
\{P_{1}, \ldots, P_{t} \}$, and for  $1 \leq i \leq t$ let $\mm_i$ be
the ideal sheaf of $P_{i}$ on $X$. Let $Z_{i}$ be the closed
subscheme of $X$ corresponding to the ideal sheaf $\mm_i^{2} \subset
\OO_X$, and define $Z := \cup_{i=1}^{t} Z_{i}$.  Note that
$H^{0}(Z,\OO_Z) \simeq \prod_{i=1}^{t} H^{0}(Z_{i},\OO_{Z_{i}}) $ and
also that $H^{0}(Z_{i},\OO_{Z_{i}}) \simeq
\OO_{X,P_{i}}/\mm_{X,P_{i}}^2$; in particular
$$
|H^{0}(Z_{i},\OO_{Z_{i}})|=
|\OO_{X,P_{i}}/\mm_{X,P_{i}}| \cdot
|\mm_{X,P_{i}}/\mm_{X,P_{i}}^2| =
p \cdot p^2 = p^{3}.
$$
Given a collection of $s$ points $W \subset X(\fp)$, we define a subset
$T(W) \subset H^0(Z, \OO_Z)$ using the above isomorphisms: let $ T(W) :=
\prod_{i=1}^t T_{i} $ where $T_{i} \subset H^{0}(Z_{i},\OO_{Z_{i}})
\simeq \OO_{X,P_{i}}/\mm_{X,P_{i}}^2$ is given by
$$
T_{i} :=
\begin{cases}
\{\phi \in \mm_{X,P_{i}}/\mm_{X,P_{i}}^2 :
\phi\ne 0 \}
&\text{if $P_i \in W$} \\
\{\phi \in \OO_{X,P_{i}}/\mm_{X,P_{i}}^2 :
\phi \not \equiv 0 \mod \mm_{X,P_{i}} \}
&\text{if $P_i  \in X(\fp) - W$}
\end{cases}
$$
Note that $|T_{i}| = p^{2}-1$ if $P_i \in W$, and $|T_{i}| =
p^{3}-p^{2}$ if $P_i \not \in W$, hence $|T(W)|=
(p^{2}-1)^{s}(p^{3}-p^{2})^{t-s}$.
We also note that if $f \in S(d)$ and $f|_Z \in T$, then $(H_{f} \cap
X)(\fp) = W$, and $H_{f}
\cap X$ is smooth at all $P_{i} \in W$.

Letting $W$ range over all subsets of $X(\fp)$ of cardinality $s$,
define a subset $T \subset H^0(Z,\OO_Z)$ by
$$
T := \cup_{W \subset X(\fp) : |W| = s} T(W).
$$
Then $|T| = \binom{t}{s} \cdot
(p^{2}-1)^{s}(p^{3}-p^{2})^{t-s}$, and, by construction,
$f|_Z \in T$ is equivalent to $X \cap H_f$ being smooth at all points
in $X(\fp)$ and that $|X \cap
H_f|=s$.

By Theorem~\ref{thm:bertini-with-taylor} and
\eqref{eq:number-of-smooth-guys}, together with
$$
\zeta_X(s) = \zeta_U(s) \prod_{P \in X(\fp)} (1-p^{-s})^{-1}
= \zeta_U(s) (1-p^{-s})^{-t}
$$
we find that for $d \to \infty$,
$$
\frac{|\{f \in S(d) : |X(\fp) \cap H_f|=s \text{ and $X \cap H_f$ is
    smooth }\}|}
{ |\{ f \in S(d) : \text{$X \cap H_f$ is    smooth }\}|}
$$
$$
=
\frac{|S(d)|/\zeta_U(3) \cdot |T|/|H^0(Z,\OO(Z))| }{|S(d)|/\zeta_X(3)}
\cdot (1+o(1))
$$
$$
=
\frac{|T|/|H^0(Z,\OO(Z))| }{(1-p^{-3})^{t}}
\cdot (1+o(1)),
$$
which, since $|T| = \binom{t}{s} \cdot
(p^{2}-1)^{s}(p^{3}-p^{2})^{t-s}$ and $|H^0(Z,\OO(Z))| = p^{3t}$,
equals
$$
\frac{\binom{t}{s} (p^{2}-1)^{s}(p^{3}-p^{2})^{t-s}}
{(p^3-1)^t}
\cdot (1+o(1))
$$
$$
=
\binom{t}{s}
\left( \frac{p+1}{p^2+p+1}\right)^s
\left( \frac{p^2}{p^2+p+1}\right)^{t-s}
\cdot (1+o(1))
$$
and hence concluding the proof of Proposition~\ref{prop:main-technical}.

\end{document}